\documentclass[11pt]{amsart}
\usepackage[english]{babel}
\usepackage[latin1]{inputenc}
\usepackage{amssymb,geometry,color,xcolor,graphicx}
\usepackage{amsmath,amsthm}
\usepackage{booktabs}
\usepackage{caption}
\usepackage{datetime}
\usepackage{dsfont}
\usepackage{amstext}
\usepackage{graphicx}
\usepackage{fancyhdr}
\usepackage{latexsym}
\usepackage{mathrsfs}
\usepackage{cases}
\usepackage{mathtools}
\usepackage{bm}
\usepackage{subcaption}
\usepackage{wrapfig}
\usepackage{tikz}\usetikzlibrary{positioning, shapes.geometric, shapes, arrows.meta, calc}
\usepackage{pgfplots}
\usepackage{ifthen}
\pgfplotsset{width=10cm,compat=1.9}

\usepackage{subcaption}

\usepackage{graphicx}
\usepackage{xcolor,listings}
\usepackage[numbered]{matlab-prettifier}
\usepackage{hyperref}
\hypersetup{
    colorlinks=true, 
    linktoc=all,     
    linkcolor=blue,  
}
\usepackage{listofitems} 
\usetikzlibrary{arrows.meta} 
\usepackage[outline]{contour} 
\contourlength{1.4pt}

\tikzset{>=latex} 

\usepackage{xcolor}
\colorlet{myred}{red!80!black}
\colorlet{myblue}{blue!80!black}
\colorlet{mygreen}{green!60!black}
\colorlet{myorange}{orange!70!red!60!black}
\colorlet{mydarkred}{red!30!black}
\colorlet{mydarkblue}{blue!40!black}
\colorlet{mydarkgreen}{green!30!black}

\tikzset{
  >=latex, 
  node/.style={thick,circle,draw=myblue,minimum size=22,inner sep=0.5,outer sep=0.6},
  node in/.style={node,green!20!black,draw=mygreen!30!black,fill=mygreen!25},
  node hidden/.style={node,blue!20!black,draw=myblue!30!black,fill=myblue!20},
  node convol/.style={node,orange!20!black,draw=myorange!30!black,fill=myorange!20},
  node out/.style={node,red!20!black,draw=myred!30!black,fill=myred!20},
  connect/.style={thick,mydarkblue}, 
  connect arrow/.style={-{Latex[length=4,width=3.5]},thick,mydarkblue,shorten <=0.5,shorten >=1},
  node 1/.style={node in}, 
  node 2/.style={node hidden},
  node 3/.style={node out}
}
\def\nstyle{int(\lay<\Nnodlen?min(2,\lay):3)} 
\def\n2style{int(\lay<\Nnodlen+4?min(2,\lay-4):3)} 

\numberwithin{equation}{section}

\newcommand{\udef}{\mathrel{\mathop:}=}

\newcommand{\R}{\mathbb{R}}

\newcommand{\N}{\mathbb{N}}

\newcommand{\de}{\,\mathrm{d}}
\usepackage{dsfont}

\theoremstyle{plain}
\newtheorem{thm}{Theorem}[section]

\newtheorem{exm}[thm]{Example}

\renewcommand{\ss}{\scriptstyle}
\def\svdots{\vbox{\baselineskip=1.5pt\lineskiplimit=0pt
	\kern1.5pt \hbox{$\ss .$}\hbox{$\ss .$}\hbox{$\ss .$}}}

\begin{document}

\title{Physics Informed Neural Networks for an Inverse Problem in Peridynamic Models
}
\author[Difonzo]{Fabio V. Difonzo}
\address{Dipartimento di Matematica, Universit\`a degli Studi di Bari Aldo Moro, Via E. Orabona 4, 70125 Bari, Italy}
\email{fabio.difonzo@uniba.it}
\author[Lopez]{Luciano Lopez}
\address{Dipartimento di Matematica, Universit\`a degli Studi di Bari Aldo Moro, Via E. Orabona 4, 70125 Bari, Italy}
\email{luciano.lopez@uniba.it}
\author[Pellegrino]{Sabrina F. Pellegrino}
\address{Dipartimento di Ingegneria Elettrica e dell'Informazione, Politecnico di Bari, Via E. Orabona 4, 70125 Bari, Italy}
\email{sabrinafrancesca.pellegrino@poliba.it}

\subjclass{34A36, 15B99}

\keywords{Physics Informed Neural Network, Peridynamic Theory, Radial Basis Functions, Inverse Problems}

\null\hfill Version of \today $, \,\,\,$ \xxivtime

\begin{abstract}
Deep learning is a powerful tool for solving data driven differential problems and has come out to have successful applications in solving direct and inverse problems described by PDEs, even in presence of integral terms. In this paper, we propose to apply radial basis functions (RBFs) as activation functions in suitably designed Physics Informed Neural Networks (PINNs) to solve the inverse problem of computing the peridynamic kernel in the nonlocal formulation of classical wave equation, resulting in what we call RBF-iPINN. We show that the selection of an RBF is necessary to achieve meaningful solutions, that agree with the physical expectations carried by the data. We support our results with numerical examples and experiments, comparing the solution obtained with the proposed RBF-iPINN to the exact solutions.
\end{abstract}

\maketitle

\pagestyle{myheadings}
\thispagestyle{plain}
\markboth{F.V. DIFONZO, L. LOPEZ AND S.F. PELLEGRINO}{RBF-iPINNs FOR AN INVERSE PROBLEM IN PERIDYNAMIC MODELS}

\section{Introduction to the peridynamic inverse problem}

We consider the following PDE in peridynamic formulation:
\begin{equation}\label{eq:periPDE}
\partial_{tt}\theta(x,t)=\int_{\R}C(|x-y|)[\theta(x,t)-\theta(y,t)]\,\de y,
\end{equation}
where $C:\R\to\R$, representing the so-called kernel function, is a nonnegative even function. \\

In the one-dimensional case, the model describes the dynamic response of an infinite bar composed of a linear microelastic material.

The main important aspect of such constitutive model is that it takes into account long-range interactions and their effects. The equation of motion~\eqref{eq:periPDE} was proposed by Silling in 2000 in~\cite{SILLING2000} in the framework of continuum mechanics theory with the name of linear peridynamic. This is an integral-type nonlocal model involving only the displacement field and not its gradient. This leads to a theory able to incorporate cracks, fractures and other kind of singularities.

The general initial-value problem associated with~\eqref{eq:periPDE} is well-posed (see~\cite{Emmrich_Puhst_2013,Emmrich_Puhst_2015}) and due to the presence of long-range forces, the solution shows a dispersive behavior. The length-scale of the long-range interactions is parameterized by a positive scalar value $\delta>0$ called horizon, which represents the maximum interaction distance between two material particles.

If the kernel function $C$, also known as micromodulus function, is a suitable generalized function, in the limit of short-range forces, or equivalently taking the limit as $\delta\to0$, the linear peridynamic model~\eqref{eq:periPDE} reduces to the wave equation $\partial_{tt}\theta(x,t)-\partial_{xx}\theta(x,t)=0$, (see~\cite{WECKNER2005,OTERKUS2014}). As a consequence, the length-scale parameter $\delta$ can be viewed as a measure of the degree of nonlocality of the model.

In order to maintain the consistency with Newton's third law, the micromodulus function must be even:
\begin{equation}
\label{eq:Ceven}
C\left(\xi\right)=C\left(-\xi\right),\quad \xi\in\R.
\end{equation}

Moreover, due to the dispersive effects $C$ must be such that
\begin{equation}
\label{eq:integral_assump}
\int_{\R} \left(1-\cos(k\xi)\right) C\left(\xi\right) \de\xi > 0,
\end{equation}
for every wave number $k\ne 0$.

Additionally, since the interaction between two material particle should become negligible as the distance between particles become very large, we can assume that
\begin{equation}
\label{eq:Climit}
\lim_{\xi\to\pm\infty} C\left(\xi\right) =0.
\end{equation}

If a material is characterized by a finite horizon, so that no interactions happen within particles that have relative distance greater than $\delta$, then we can assume that the support of the kernel function is given by $[-\delta,\delta]$ and the model~\eqref{eq:periPDE} writes as
\begin{equation}
\label{eq:periPDEsupport}
\partial_{tt}\theta(x,t)=\int_{B_{\delta}(x)}C(\|x-y\|)[\theta(x,t)-\theta(y,t)]\,\de y.
\end{equation}
Of course, such condition is less restrictive than~\eqref{eq:Climit}.

It is clear that a different microelastic material corresponds to a different kernel function and, as a consequence, the kernel function involved in the model provides different constitutive models.

In this paper, we aim to solve the inverse problem described in \eqref{eq:periPDE} for learning the shape of the kernel function $C$, by implementing
a Physics Informed Neural Network (PINN). More specifically, we show that this inverse problem requires a careful selection of activation functions in all the layers and a correct interaction with kernel initializers. It can be seen, in fact, that a naive choice on these functions would result in unreliable predictions and possibly unfeasible solutions. More precisely, we see that, if the neural network structure is not chosen accordingly to appropriate geometric knowledge relative to the data, then PINN output returns different, still acceptable, results, showing a lack of uniqueness. Therefore we will show that, as soon as the peridynamic operator is bounded on a compact support $[-\delta,\delta]$ and the PINN architecture is build accordingly, as a consequence of the well posedness conditions of the peridynamic formulation, the learned kernel fulfills all the requirements expected, provided that PINN structure is accurate enough.

\section{Introduction to PINNs}

Physics-Informed Neural Networks (PINNs) have emerged as a transformative approach to tackle both direct and inverse problems associated with PDEs. These innovative neural network architectures seamlessly integrate the principles of physics into the machine learning framework. By doing so, PINNs offer a promising solution to efficiently and accurately model, simulate, and optimize complex systems governed by PDEs. More specifically, they can be employed to solve both direct and inverse problems; in the latter case, such PINNs are commonly referred to as inverse PINNs. \\

Direct problems involve finding solutions to PDEs that describe the evolution of physical systems under specified initial and boundary conditions. Traditional numerical methods, such as finite element analysis (see~\cite{Galvanetto2018,KilicMadenci}), finite difference methods with composite quadrature formulas (see~\cite{CFLMP,Pellegrino2020}) and applied to spectral fractional models (see ~\cite{DifonzoGarrappa2023}), model order reduction methods (see~\cite{REGAZZONI2019}), meshfree methods (see~\cite{Silling_Askari_2005,Galvanetto2016}), adaptive refinement techniques (see~\cite{Bobaru_2009,ALEBRAHIM2023109710}) and collocation and Galerkin methods (see~\cite{ALEBRAHIM2023}) have been widely used for solving direct problems. Moreover, more recently spectral methods with volume penalization techniques (see for instance~\cite{LP,LP2021,Du,Jafarzadeh}) and Chebyshev spectral methods (see~\cite{LPeigenv,LPcheby,LPcheby2022,Jafarzadeh2023}) have been developed in order to increase the order of convergence, to improve the accuracy of the results and to maintain the consistency of the method even in presence of singularities.

However, these approaches often require substantial computational resources and may struggle with high-dimensional or non-linear problems. Additionally, such methods need to know the constitutive parameters of the model to predict fractures in the material under consideration and, in suitable configurations, they fail to impose boundary conditions. In order to provide some hint in this direction, a data-driven approach can be developed. In~\cite{SUKUMAR2022} the authors propose a geometry-aware method in physics informed neural network to exactly imposing boundary conditions over complex domains. In~\cite{Zhou2023} the authors investigate both a forward and an inverse problems of high-dimensional nonlinear wave equations via a deep neural networks with the activation function. While in~\cite{Zunino} a combination of an orthogonal decomposition with a neural network is applied to build a reduced order model. In~\cite{MadenciCMAME2023}, the authors present an unsupervised convolutional neural network architecture with nonlocal interactions for solving PDEs using Peridynamic Differential Operator as a convolutional filter. Indeed, this approach results to be very efficient when the model is governed by an integral operator (see for instance~\cite{YuanJCP2022}).

Inverse problems, on the other hand, are concerned with determining unknown parameters, boundary conditions, or the PDE itself, given limited or noisy observations of the system behavior. These problems frequently arise in real-world applications, including medical imaging \cite{ChenEtAl20}, geophysics \cite{bandaiGhezzehei2022,depinaEtAl2022,BERARDI2024}, material characterization \cite{XU2023}, and industrial process optimization \cite{MENG2023}. Inverse problems are inherently ill-posed, as multiple solutions or no solutions may exist, making their resolution challenging. In fact, several issues could arise in solving inverse problems, especially related to irregular geometries \cite{GAO2022}, or also small data regimes, incomplete data or incomplete models \cite{RAISSI2019}.

In the context of nonlocal elasticity theory, in~\cite{Turner} the authors propose a methodology based on a constrained least squares optimization to solve inverse problems in heterogeneous media using state-based peridynamics in order to derive parameter values characterizing several material properties and to establish conditions for fracture patterns in geological setting.

Thus, Physics-Informed Neural Networks represent a paradigm shift in the way to approach direct and inverse problems associated with PDEs. Their ability to combine data-driven learning with physical principles opens up new frontiers in scientific discovery, engineering design, and problem-solving across a wide spectrum of domains.

\subsection{PINN paradigm}

In this paper, we will consider a Feed-Forward fully connected Neural Network (FF-DNN), also called Multi-Layer Perceptron (MLP) (see \cite{bengio2003,Cybenko1989} and references therein). \\
In a PINN the solution space is approximated through a combination of activation functions, acting on all the hidden layers, with the independent variable used as the network inputs. Letting $x\in\R^n$, in a Feed-Forward network each layer feeds the next one through nested transformation, so that a it can be seen, letting $L$ be the number of layers, as
\begin{equation}\label{eq:NN}
\begin{aligned}
z_0 &= x, \\
z_l &= \sigma_l\left(\Lambda_l(z_{l-1})\right),\,\,\Lambda_l(z_{l-1})\udef W_lz_{l-1}+b_l,\quad l=1,\ldots,L,
\end{aligned}
\end{equation}
where, for each layer $l=1,\ldots,L$, $\sigma_l:\R^n\to\R^m$ is the activation function, which operates componentwise, $W_l$ is the weight matrix and $b_l$ is the bias vector. Thus, the output $z_L\in\R^m$ of a FF-NN can be expressed as a single function of the input vector $x$, defined as the composition of all the layers above in the following way:
\[
z_L\udef(\sigma_L\circ\Lambda_L\circ\ldots\circ\sigma_1\circ\Lambda_1)(x).
\]
The aim of a PINN is to minimize, through a Stochastic Gradient Descent method, a suitable objective function called \emph{loss function}, that would take into account the physics of the problem, with respect to all the components, called trainable parameters, of $W_l,b_l$, for $l=1,\ldots,L$.

More specifically, given a general PDE of the form $\mathcal{P}(f)=0$, where $\mathcal{P}$ represents the differential operator acting on $f$, the loss function used by a PINN is usually given by
\begin{equation}\label{eq:lossPINN}
\mathcal{L}(f)\udef\|f-f^*\|+\|\mathcal{P}(f)-0^*\|,
\end{equation}
where $f^*$ is the training dataset (of points inside the domain or on the boundary), and $0^*$ is the expected (true) value for the differential operation $\mathcal{P}(f)$ at any given training or sampling point; the chosen norm $\|\cdot\|$ (it may be different for each term in the loss function) depends on the functional space and the specific problem. Selecting a correct norm (so to avoid overfitting) for the loss function evaluation is an important problem in PINN, and recently in \cite{TAYLOR2023} authors have proposed spectral techniques based on Fourier residual method to overcome computational and accuracy issues. The first term in the right-hand side of \eqref{eq:lossPINN} is referred to as data fitting loss, while the second term is referred to as residual loss, which is responsible to make a NN be informed by physics.

The operator $\mathcal{P}$ is usually performed using \texttt{autodiff} (Automatic Differentiation algorithm). In the context of peridynamic theory, in \cite{HAGHIGHAT2021} authors propose, for the first time, a nonlocal alternative to \texttt{autodiff} by replacing the evaluation of $f$ and its partial derivatives through the action of a Peridynamic Differential Operator (PDDO) on $f$.

A recent review on PINNs and related theory can be found in \cite{Cuomo2022}.

\section{RBF-iPINN for the kernel function}

In case one wants to solve an inverse problem then there will be more trainable parameters than only those coming from weight matrices and bias vectors. Hence, such further parameters have to be considered in the minimization iterations and their respective gradients must be computed as well. \\
However, in our case, the inverse problem does not involve the mere computation of scalar quantities, but rather a whole function, specifically the kernel function $C$ in \eqref{eq:periPDE}, which has also analytical and geometrical properties to be accounted for, such as nonnegativity and symmetry. In the PINN architecture proposed, these features reflect in the implementation of a NN model with two separated sets of layers, one for $C$ and the other for $\theta$, and whose output would be both the solution to \eqref{eq:periPDE} and the unknown function $C$; moreover, the loss function \eqref{eq:lossPINN} has been accordingly endowed with further terms necessary to enforce the requirements on $C$.

From the point of view of the architecture, while nonnegatitivity of $C$ has been simply enforced by requiring all trainable parameters in \eqref{eq:NN} to be nonnegative, symmetry has required a more specific treatment, both in terms of activation functions and in the way we have defined the loss function. Our idea has been to wisely select different activation functions for the two different sets of layers, inspired by the properties coming along with $C$ and the data on $\theta$. In the following section we will introduce and discuss the technical approaches to deal with activation and loss functions.

\subsection{Radial Basis Function Layer}

As activation function for the first layer, whose input is $x$, a Radial Basis Function (RBF) is selected. By definition, a Radial Basis Function (RBF) is a real-valued function whose output depends only on the distance from a fixed center or prototype point. An RBF can be defined as

\begin{equation}\label{eq:rbf_class}
\phi(x) = \phi(\|x - c\|),
\end{equation}

where $\phi$ is the RBF function, $x$ is the input to the RBF, $c$ is the center or prototype point,  $\|x - c\|$ represents the distance between $x$ and $c$.

RBFs are commonly used in various fields and, when used as activation functions in neural networks, they give rise to Radial Basis Function Neural Networks (RBFNNs) (see \cite{fasshauer2007meshfree}). \\

We have experimented on two families of RBFs \eqref{eq:rbf_class}, given by
\begin{subequations}\label{eq:rbf}
\begin{align}
\sigma_\textrm{rbf}(x) &\udef \frac{\rho}{1+\gamma(x-\mu)^2},\quad\rho,\gamma,\mu>0, \label{eq:rbfGauss} \\
\sigma_\textrm{rbf}(x) &\udef \rho\sqrt{1+\gamma(x-\mu)^2},\quad\rho,\gamma,\mu>0,  \label{eq:rbfV}
\end{align}
\end{subequations}
called inverse quadratic and  multiquadric RBFs respectively, where all the parameters above could be taken to be trainable. This approach has recently been proposed in \cite{BaiEtAl2023} in the context of direct problem for nonlinear PDEs, but used in the middle layer. However, our experiments solving the inverse problem of learning the kernel function of a PDE in the peridynamic formulation, we have extensively noticed that such a choice is not efficient, resulting in nonphysical results and excessively large computational time and cost. In fact, this has led us to introduce an RBF inverse PINN, that we called RBF-iPINN, where the first layer is implemented with an RBF activation function, that has significantly sped up performance while providing the expected result if compared to the exact solution.

\begin{figure}
\centering
\begin{subfigure}{.32\textwidth}
\begin{tikzpicture}[
  declare function={
    func(\x)= 1/(1+\x^2);
  }
]
\begin{axis}[
    width = \linewidth,
    axis lines = left,
    xlabel = \(x\),
    ylabel = {\(\sigma_{\mathrm{rbf}}(x)\)},
    ymin = 0,
]

\addplot[
    domain=-10:10,
    samples=200,
    color=black,
] {func(x)};
\end{axis}
\end{tikzpicture}
\end{subfigure}
\qquad\qquad
\begin{subfigure}{.32\textwidth}
\begin{tikzpicture}[
  declare function={
    func(\x)= sqrt(1+x^2);
  }
]
\begin{axis}[
    width = \linewidth,
    axis lines = left,
    xlabel = \(x\),
    ylabel = {\(\sigma_{\mathrm{rbf}}(x)\)},
    ymin = 0,
]

\addplot[
    domain=-10:10,
    samples=200,
    color=black,
] {func(x)};
\end{axis}
\end{tikzpicture}
\end{subfigure}
\caption{Qualitative shapes of Radial Basis Functions defined in \eqref{eq:rbfGauss} and \eqref{eq:rbfV}, respectively.}
\label{fig:influenceFunction}
\end{figure}
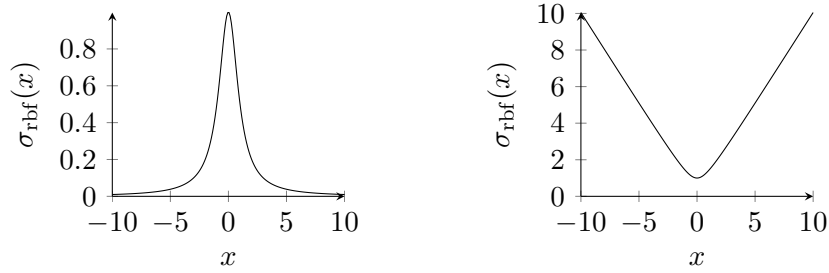

\subsection{RBF-iPINN Structure}

Since the kernel function $C$ is a variable of the sole spatial variable $x$, while the solution to \eqref{eq:periPDE} $\theta$ depends on both space and time, then $x,t$ need to be handled separately. To this purpose, the proposed RBF-iPINN is implemented in a serialized fashion That is, its first layer is activated by an RBF as in \eqref{eq:rbf} with the only input $x$, followed by 8 layers with 20 neurons each, activated by \texttt{ReLu} function. The output of this sequence of layer is then concatenated with $t$, to produce the input for the second portion on the Neural Network; more specifically, we put again 8 layers with 20 neurons each, having \texttt{sigmoid} as their activation function, before providing the overall output of the RBF-iPINN, which is thus represented by a tensor for $C$ and a tensor for $\theta$. The structure is sketched in Figure \ref{fig:NNstructure}. Moreover, the two sets of full hidden layers are endowed with a kernel initializer of type \texttt{glorot\_normal} and \texttt{random\_uniform} respectively; also, a kernel regularizer of type \texttt{l1\_l2}, with weights $l_1=l_2=0.01$ is applied to the first portion of the NN, which yields the computed kernel function.

\subsection{Loss Function}

Given the constraints on $C$, we have to accordingly construct a loss function as in \eqref{eq:lossPINN} with as many components as there are constraints to be enforced in the RBF-PINN. More specifically, we consider the following components to be part of $\mathcal{L}$:
\begin{subequations}\label{eq:loss}
\begin{align}
    &\textrm{pde loss: }\mathcal{L}_\textrm{pde}\udef\|\mathcal{P}(f)-0^*\|_2, \label{eq:loss_pde}\\
    &\textrm{data fitting loss: }\mathcal{L}_\textrm{data}\udef\|f-f^*\|_2, \label{eq:loss_data}\\
    &\textrm{symmetry loss: }\mathcal{L}_\textrm{sym}\udef\|f(x,t)-f(-x,t)\|_1. \label{eq:loss_sym}
\end{align}
\end{subequations}
In the selection of norms above, we have been guided by the features we want our PINN to take into account. More specifically, since we are interested in fitting data and in satisfying our model as much as possible, we selected the $2$-norm for these contributions; however, since symmetry is to be expected from the PINN architecture and, in particular, from the first RBF layer, we measure its loss via the $1$-norm, which is more sensible to small errors. Finally, we consider a weighted sum of the contributions given above as
\begin{equation}\label{eq:PINNLoss}
\mathcal{L}=w_\textrm{pde}\mathcal{L}_\textrm{pde}+w_\textrm{data}\mathcal{L}_\textrm{data}+w_\textrm{sym}\mathcal{L}_\textrm{sym}.
\end{equation}

\subsection{Learning rate}\label{subsec:lr}

The learning rate $\alpha$ has been selected to be decreasing with the epoch in a quadratic way. More precisely, we implemented the following scheduler:
\begin{equation}\label{eq:lr}
\begin{aligned}
\alpha_0 &= 10^{-4}, \alpha_1 = 0.7\alpha_0 \\
\alpha_i &= \left(1-\left(\frac{i}{N}\right)^2\right)\alpha_0+\left(\frac{i}{N}\right)^2\alpha_1,\quad i=1,\ldots,N,
\end{aligned}
\end{equation}
where $N$ is the number of epochs chosen for the training.

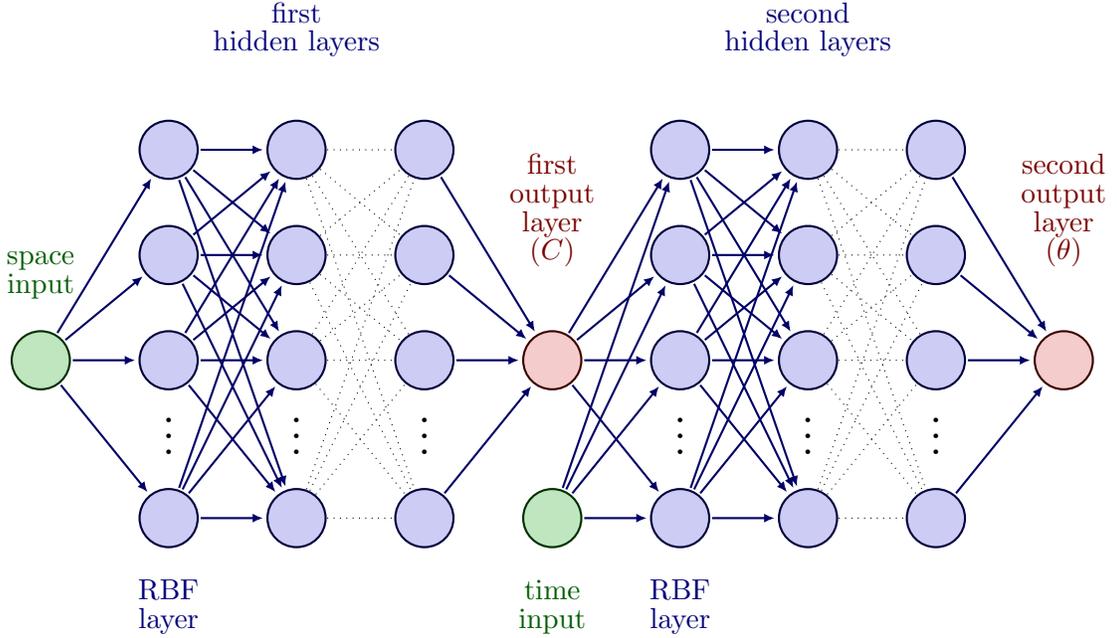
\begin{figure}
\begin{tikzpicture}[x=1.7cm,y=1.4cm]
  \message{^^JNeural network, shifted}
  \readlist\Nnod{1,4,4,4,1} 
  \def\yshift{0.5} 

  \message{^^J  Layer}
  \foreachitem \N \in \Nnod{ 
    \def\lay{\Ncnt} 
    \pgfmathsetmacro\prev{int(\Ncnt-1)} 
    \message{\lay,}
    \foreach \i [evaluate={\c=int(\i==\N); \y=\N/2-\i-\c*\yshift;
                 \x=\lay; \n=\nstyle;}] in {1,...,\N}{ 
      \node[node \n] (N\lay-\i) at (\x,\y) {}
      ;
      \ifnum\lay>1 
      \ifnum\lay<4
        \foreach \j in {1,...,\Nnod[\prev]}{ 
          \draw[connect arrow] (N\prev-\j) -- (N\lay-\i);
        }
      \fi
      \fi 
      \ifnum\lay>4
        \foreach \j in {1,...,\Nnod[\prev]}{ 
          \draw[connect arrow] (N\prev-\j) -- (N\lay-\i);
        }
      \fi
      \ifnum\lay=4
        \foreach \j in {1,...,\Nnod[\prev]}{ 
          \draw[dotted] (N\prev-\j) -- (N\lay-\i);
        }
      \fi

    }
    \ifnum\lay>1 \ifnum\lay<5
    \path (N\lay-\N) --++ (0,1+\yshift) node[midway,scale=1.5] {$\vdots$};
    \fi\fi
  }

  \node[above=0.2,align=center,mygreen!60!black] at (N1-1.90) {space\\[-0.2em]input};
  \node[below=0.2,align=center,myblue!60!black] at (N2-4.-90) {RBF\\[-0.2em] layer};
  \node[above=0.5,align=center,myblue!60!black] at (N3-1.90) {first\\[-0.2em]hidden layers};
  \node[above=0.5,align=center,myred!60!black] at (N\Nnodlen-1.90) {first\\[-0.2em]output\\[-0.2em]layer\\[-0.2em]($C$)};

  \node[node 1] (N5-2) at (5,-2.5) {};
  \node[below=0.2,align=center,mygreen!60!black] at (N5-2.-90) {time\\[-0.2em]input};

  \readlist\Nnod{2,4,4,4,1} 

  \message{^^J  Layer}
  \foreachitem \N \in \Nnod{ 
    \pgfmathsetmacro\lay{int(\Ncnt+4)} 
    \pgfmathsetmacro\prev{int(\Ncnt+3)} 
    \pgfmathsetmacro\prevv{int(\Ncnt-1)}
    \message{\lay,}
    \foreach \i [evaluate={\c=int(\i==\N); \y=\N/2-\i-\c*\yshift;
                 \x=\lay; \n=\n2style;}] in {1,...,\N}{ 
      \ifnum \lay>5
      \node[node \n] (N\lay-\i) at (\x,\y) {};
      \fi
      \ifnum\lay>5 
      \ifnum\lay<8
        \foreach \j in {1,...,\Nnod[\prevv]}{ 
          \draw[connect arrow] (N\prev-\j) -- (N\lay-\i);
        }
      \fi
      \fi 
      \ifnum\lay>8
        \foreach \j in {1,...,\Nnod[\prevv]}{ 
          \draw[connect arrow] (N\prev-\j) -- (N\lay-\i);
        }
      \fi
      \ifnum\lay=8
        \foreach \j in {1,...,\Nnod[\prevv]}{ 
          \draw[dotted] (N\prev-\j) -- (N\lay-\i);
        }
      \fi

    }
    \ifnum\lay>5 \ifnum\lay<9
    \path (N\lay-\N) --++ (0,1+\yshift) node[midway,scale=1.5] {$\vdots$};
    \fi\fi
  }

  \node[below=0.2,align=center,myblue!60!black] at (N6-4.-90) {RBF\\[-0.2em] layer};
  \node[above=0.5,align=center,myblue!60!black] at (N7-1.90) {second\\[-0.2em]hidden layers};
  \node[above=0.5,align=center,myred!60!black] at (N9-1.90) {second\\[-0.2em]output\\[-0.2em]layer\\[-0.2em]($\theta$)};

\end{tikzpicture}
\caption{RBF-iPINN structure. In our simulations, we set the number of hidden layers, after each RBF layer, to $8$, and the number of neurons per layer to $20$.}
\label{fig:NNstructure}
\end{figure}

\section{Numerical Simulations}

In this section, we show results with our RBF-PINN. All the experiments have been run using 1000 epochs with a learning rate defined in \eqref{eq:lr} and employed the ADAM optimizer. The machine used for the experiments is an Intel Core i7-8850H CPU at 2.60GHz and 64 GB of RAM. Moreover, the PINNs, providing results in the examples below, have been developed in the python library TensorFlow using the interface Keras. \\
Moreover, real data, which are used in our simulations to compute the loss function \eqref{eq:loss_data}, are synthetically built using appropriate spectral methods \cite{LP} to solve \eqref{eq:periPDE}. \\

A main feature of the numerical computation is the evaluation of the integral on the right-hand side of \eqref{eq:periPDE}. In fact, in order to exploit the power of Keras, we notice that
\begin{align*}
\int_\R C(\|x-y\|[\theta(x,t)-\theta(y,t)]\de y &= \theta(x,t)\int_\R C(\|x-y\|)\de y-\int_\R C(\|x-y\|)\theta(y,t)\de y \\
&= \theta(x,t)\int_\R C(\|x-y\|)\de y-C(\|x\|)\ast\theta(x,t),
\end{align*}
where the second term in the right-hand side above is the convolution product between the kernel $C$ and the unknown function $\theta$. It has to be noticed here that the kernel function $C$ is compactly supported, with support $[-\delta,\delta]$. Now, in order to numerically compute such convolution product, let $[0,X]$ be the space interval and let $0<x_1<x_2<\ldots<x_{N-1}<x_N=X$ be the uniform spatial discretization of the interval $[0,X]$ with stepsize $h>0$. the convolution product above can be numerically treated by determining the exact number of components in the vector $[C(x_i)]_{i=1}^n$ so that only points $x_i$ such that
\[
|x_i-x_j|<\delta,\,\,i,j=1,\ldots,N,
\]
come into play when computing $C(\|x\|)\ast\theta(x,t)$. Since $x_i=i\cdot h$, then we deduce that the only indices involved in the convolution product are $i,j=1,\ldots,N$ such that
\[
|i-j|<\frac{\delta}{h}.
\]

Since the peridynamic integral-operator in~\eqref{eq:periPDEsupport} is linear, we can derive in terms of the Green's function a solution of the initial-value problem using continuous Fourier Transform. Such solution can be used to provide a dataset for the next simulations.

In the next two examples, we exemplify on datasets derived from V-shaped kernel functions. This choice of kernel is justified by the fact that it is implemented in some nonlocal formulations of Richards' equation as it is able to easily incorporate Dirichlet boundary condition in the model (see~\cite{BDFP}).

Thus, we select as activation function for the first layer of our PINN model an RBF of type \eqref{eq:rbfV}. Further, we tried both to keep all the three parameters $\gamma,\rho,\mu$ trainable and to fix $\gamma$ while letting $\rho,\mu$ be trainable. According to our experience and for the following two cases, fixing $\gamma$ improves convergence performance and result quality.

\begin{exm}\label{subsec:ex1}

Here we consider a dataset with $t\in[0,20], x\in[-10,10]$ with spatial stepsize $h=2\cdot10^{-1}$ and $\delta=10$, and for which the analytical expression of the kernel is
\begin{equation}\label{eq:kernelType2}
C(x)=\frac35|x|.
\end{equation}
We set $\gamma=0.09$, obtaining the results showed in Figure \ref{fig:learnedKernel_type2_gamma009}, where we compare the true kernel function in \eqref{eq:kernelType2} to the output of the proposed inverse RBF-iPINN.
Setting $\gamma=0.05$ in \eqref{eq:rbfV} provides qualitatively comparable results, as can be observed in Figure \ref{fig:learnedKernel_type2_gamma005}.
\begin{figure}
    \centering
     \begin{subfigure}{0.48\textwidth}
    \includegraphics[width=\textwidth]{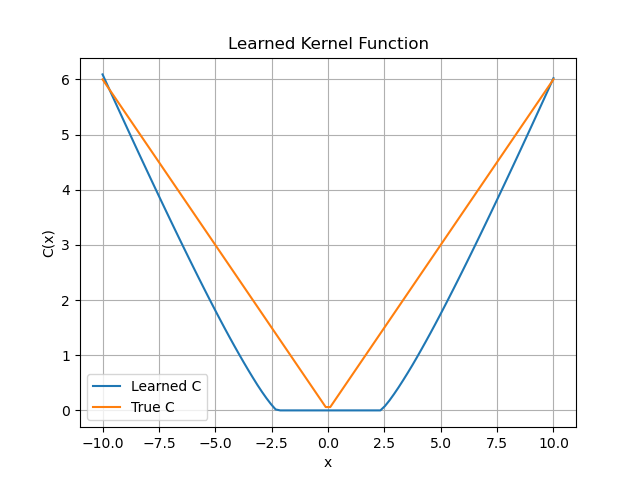}
    \caption{$\gamma=0.09$}
    \label{fig:learnedKernel_type2_gamma009}
    \end{subfigure}
     \begin{subfigure}[b]{0.48\textwidth}
         \includegraphics[width=\textwidth]{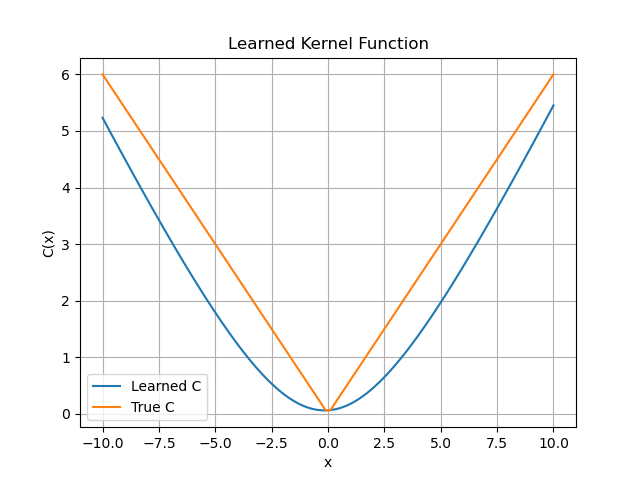}
    \caption{$\gamma=0.05$}
    \label{fig:learnedKernel_type2_gamma005}
     \end{subfigure}
     \caption{Learned kernel functions relative to Example \ref{subsec:ex1} for different values of $\gamma$ in \eqref{eq:rbfV}.}
     \label{fig:learnedKernel_type2}
\end{figure}
\end{exm}

\begin{exm}\label{subsec:ex2}

Here we consider a dataset with $t\in[0,20], x\in[-10,10]$ with spatial stepsize $h=2\cdot10^{-1}$ and $\delta=1$, and for which the analytical expression of the kernel is
\begin{equation}\label{eq:kernelType3}
C(x)=
\begin{cases}
\frac{\delta-x-10}{\delta}, & x\leq-10+\delta, \\
0, & -10+\delta<x\leq10-\delta, \\
\frac{\delta+x-10}{\delta}, & x>10-\delta.
\end{cases}
\end{equation}
We set $\gamma=0.09$, obtaining the results showed in Figure \ref{fig:learnedKernel_type3_gamma009}, where we compare the true kernel function in \eqref{eq:kernelType3} to the output of the proposed inverse RBF-iPINN. Setting $\gamma=0.05$ provides results in Figure \ref{fig:learnedKernel_type3_gamma005}.

\begin{figure}
    \centering
    \begin{subfigure}{0.48\textwidth}
         \includegraphics[width=\textwidth]{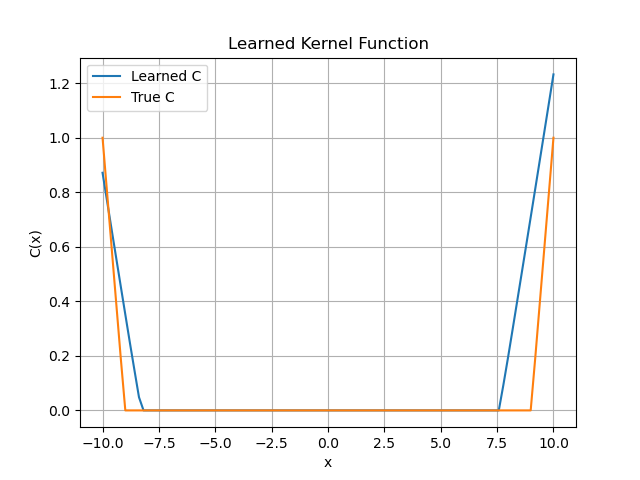}
    \caption{$\gamma=0.09$}
    \label{fig:learnedKernel_type3_gamma009}
    \end{subfigure}
    \begin{subfigure}{0.48\textwidth}
    \includegraphics[width=\textwidth]{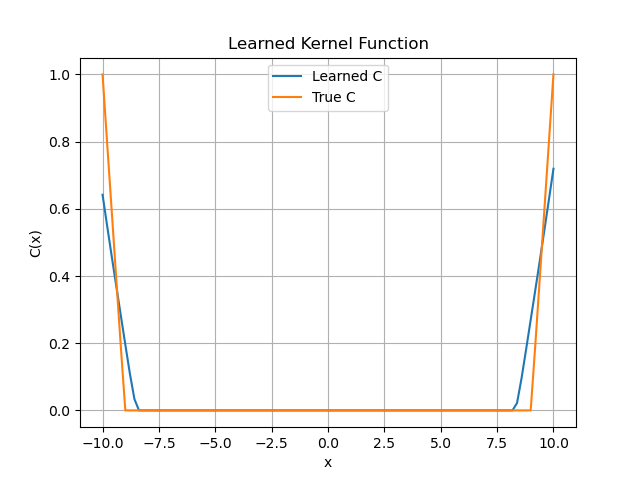}
    \caption{$\gamma=0.05$}
    \label{fig:learnedKernel_type3_gamma005}
     \end{subfigure}
     \caption{Learned kernel functions for different values of $\gamma$ in \eqref{eq:rbfV} relative to Example \ref{subsec:ex2}.}
     \label{fig:learnedKernel_type3}
\end{figure}
\end{exm}

\begin{exm}\label{subsec:ex3}

For this example we consider a bell-shaped kernel function to test the proposed RBF-iPINN. First, we tuned hyperparameters, setting the kernel regularizers \texttt{l1\_l2} with weights $0.01$ and $0.1$ respectively, in order to try and catch, as better as possible, the correct qualitative behavior of the kernel in the interior of its compact support; moreover, on the account of the knowledge of the kernel shape, we decided to activate the first layer of the RBF-iPINN through \eqref{eq:rbfGauss}, where we set the hyperparameter $\gamma=1$; finally, for a better data fitting, we also selected the sup-norm in \eqref{eq:loss_data}. \\
In this case, the neural network shows a discrete ability to catch shape and support of the bell-shaped kernel, but fails in a good approximation of the characteristic parameters, as shown in Figure \ref{fig:learnedKernel_type6}. In fact, in this case the true kernel is given by
\begin{equation}
    C(x)=\frac{4}{\sqrt{\pi}}e^{-x^2}.
\end{equation}
\begin{figure}
    \centering
    \includegraphics[width=\textwidth]{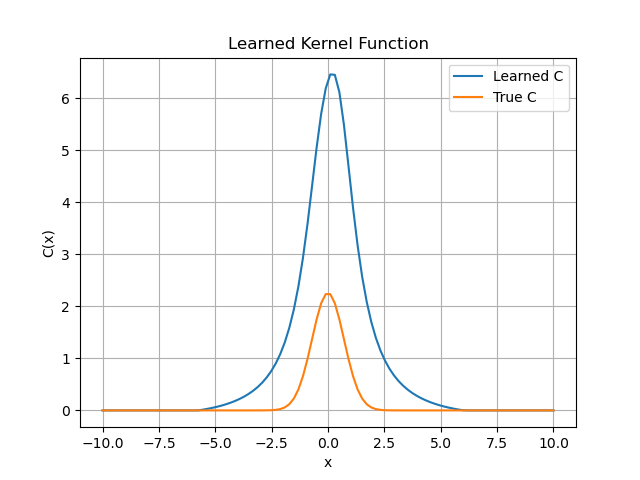}
    \caption{Learned kernel function compared to the true one from Example \ref{subsec:ex3}.}
    \label{fig:learnedKernel_type6}
\end{figure}
Therefore, we have performed a further analysis by implementing a standard inverse PINN to learn parameters $\gamma^*$ and $\sigma^*$ in
\begin{equation}\label{eq:guassianKernel}
C^*(x)\udef\gamma^*e^{-\sigma^*x^2}.
\end{equation}
Starting with initial guesses for $\gamma^*=3$ and for $\sigma^*=0.5$, we run a PINN whose structure is the same as the second portion relative to $\theta$ of the RBF-iPINN described above (see the architecture in Figure \ref{fig:NNstructure}). The training phase has been performed over $1000$ epochs and with the same learning rate scheduler described in Section \ref{subsec:lr}, but with a faster descent obtained by setting $\alpha_0=10^{-3}$. Results are depicted in Figure \ref{fig:learnedKernel_gaussian_type6}. It can be deduced that the inverse PINN has been able to correctly detect the learned values which, at convergence, are given by $\gamma^*=2.3302033$ and $\sigma^*=1.0218402$, being $\frac{4}{\sqrt{\pi}}\approx2.2567583$.

\begin{figure}
    \centering
    \includegraphics[width=\textwidth]{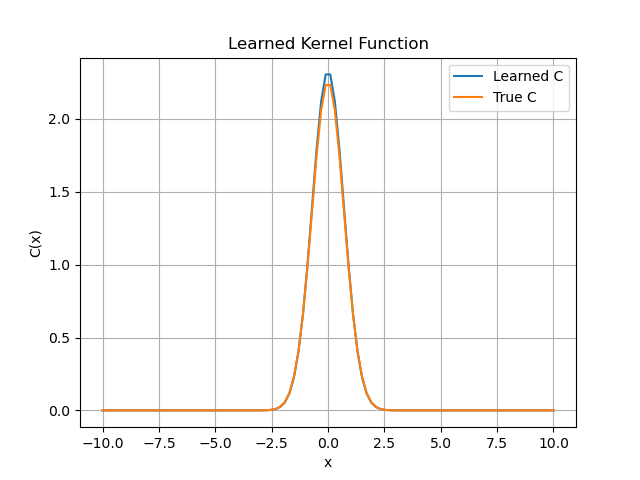}
    \caption{Gaussian kernel function learned from  \eqref{eq:guassianKernel}, compared to the true one from Example \ref{subsec:ex3}.}
    \label{fig:learnedKernel_gaussian_type6}
\end{figure}
\end{exm}

\section{Conclusions}

In this work we have analyzed a peridynamic formulation of a classical wave equation, trying to compute the kernel function responsible for the nonlocal behavior of the model considered. We have proposed to utilize a Radial Basis Function (RBF) as activation function for the first layer of a suitably designed Physics Informed Neural Network (PINN) to solve the inverse problem. Specifically, our inverse PINN architecture has two neural networks working in series: the first set of layers, whose first layer activated by some Radial Basis Function, takes the spatial data as input and yields the first output for recovering the kernel function; then, this first output is concatenated to the temporal data, thus providing a new tensor serving as input to the second set of layers, which produces the output describing $\theta$, solution to the peridynamic wave equation. We called the proposed model RBF-iPINN. We have shown that, with a wise scheduler for the learning rate and necessary initializations of the first set of layers responsible for the kernel function computation, for standard selections of the activation RBF, the RBF-iPINN can provide reliable prediction of the kernel function, in case it has a V-shape behavior. We also tackle the case of Gaussian kernel function: here, RBF-iPINN is able to adequately detect the shape, but a further analysis is necessary to learn the expected parameters of a bell-shaped function. We did so by implementing a standard inverse PINN, practically tuning the very same second set of layer of the RBF-iPINN. \\
Such models turn out to be promising tools for investigating optimal controls problems in dimension $2$ or higher (see, e.g., \cite{BerardiDifonzoGuglielmi2023}), where an inverse PINN approach seems to provide robust and scalable results. Such considerations pave the way to further investigations about how to deal with more complicated peridynamic models via PINNs and Radial Basis Functions, that seem to be a powerful approach to this kind of problems, due to their inherently symmetric nature.

\section*{Acknowledgments}

FVD has been supported by \textit{REFIN} Project, grant number 812E4967 funded by Regione Puglia.

SFP has been supported by \textit{REFIN} Project, grant number D1AB726C funded by Regione Puglia, and by \textit{PNRR MUR - M4C2} project, grant number N00000013 - CUP D93C22000430001.

The three authors gratefully acknowledge the support of INdAM-GNCS 2023 Project, grant number CUP$\_$E53C22001930001.

\bibliographystyle{plain}
\bibliography{peridynamicPinn.bib}

\end{document}